\documentclass[11pt]{amsart}

\usepackage{amsmath}
\usepackage{fullpage}
\usepackage{xspace}
\usepackage[psamsfonts]{amssymb}
\usepackage[latin1]{inputenc}
\usepackage{graphicx,color}
\usepackage[curve]{xypic} 
\usepackage{hyperref}
\usepackage{graphicx}
\usepackage[greek,english]{babel}

\usepackage{amsmath}%
\usepackage{amsthm}%
\usepackage{amscd}
\usepackage{amsfonts}%
\usepackage{amssymb}%
\usepackage{graphicx}

\usepackage{mathrsfs}

\usepackage{tikz}
\usetikzlibrary{matrix,arrows}

\usepackage{tikz-cd}

\newtheorem{theorem}{Theorem}[section]

\newtheorem{lemma}[theorem]{Lemma}

\newtheorem{proposition}[theorem]{Proposition}

\theoremstyle{remark}

\numberwithin{equation}{section}

\newcommand{\Cfrak}{\mathfrak{C}}

\newcommand{\Lfrak}{\mathfrak{L}}
\newcommand{\Ifrak}{\mathfrak{I}}
\newcommand{\Rfrak}{\mathfrak{R}}

\newcommand{\Xcal}{\mathscr{X}}

\newcommand{\Pro}{\mathbb{P}}

\newcommand{\Z}{\mathbb{Z}}

\newcommand{\F}{\mathbb{F}}
\newcommand{\Q}{\mathbb{Q}}
\newcommand{\R}{\mathbb{R}}
\newcommand{\N}{\mathbb{N}}

\newcommand{\rad}{\mathrm{rad}}

\makeatletter
\@namedef{subjclassname@2020}{%
  \textup{2020} Mathematics Subject Classification}
\makeatother

  \DeclareFontFamily{U}{wncy}{}
    \DeclareFontShape{U}{wncy}{m}{n}{<->wncyr10}{}
    \DeclareSymbolFont{mcy}{U}{wncy}{m}{n}
    \DeclareMathSymbol{\Sha}{\mathord}{mcy}{"58}

\begin{document}
\title[]{Why should $\mathbb{Z}$ be Diophantine in $\mathbb{Q}$?}

\author{Alexander Kozachinskiy}
\address{CENIA. Chile}
\email[A. Kozachinskiy]{alexander.kozachinskyi@cenia.cl}%

\author{Hector Pasten}
\address{ Departamento de Matem\'aticas,
Pontificia Universidad Cat\'olica de Chile, and CENIA. Chile}
\email[H. Pasten]{hector.pasten@uc.cl}%

\author{Xavier Vidaux}
\address{ Departamento de Matem\'aticas,
Universidad de Concepci\'on. Chile}
\email[X. Vidaux]{xvidaux@udec.cl}%


\thanks{A.K. was supported by ANID Fondecyt Iniciaci\'{o}n grant 11250060 and  National Center for Artificial Intelligence CENIA FB210017. H.P. was supported by ANID Fondecyt Regular grant 1230507 from Chile. X.V. was supported by Proyecto 2025001477INV from Universidad de Concepci\'on.}
\date{\today}
\subjclass[2020]{Primary: 11U05; Secondary: 11G05, 11N32} %
\keywords{Diophantine sets, rational numbers, Hilbert's tenth problem}%

\begin{abstract} We give a Diophantine definition of the integers in the rational numbers conditional on the $abc$ conjecture, Dickson's conjecture, and finiteness of $2$-primary parts of Shafarevich--Tate groups. We also discuss some consequences of the existence of a Diophantine definition of the integers in the rational numbers. The goal of this work is to serve as an entry point to recent developments around Hilbert's tenth problem for the rational numbers.
\end{abstract}

\maketitle



\section{Introduction}

\subsection{The result} Let us consider the following three conjectures:
\begin{itemize}
\item[(X)] Finiteness of the $2$-primary part of Shafarevich--Tate groups of elliptic curves over $\Q$.
\item[(Y)] The $abc$ conjecture.
\item[(Z)] Dickson's conjecture on simultaneous prime values of linear polynomials.
\end{itemize}

The main result of this note is:

\begin{theorem}\label{ThmMain} Assuming (X), (Y), and (Z), there is a Diophantine definition of $\Z$ in $\Q$.
\end{theorem}

As a consequence, under (X), (Y), and (Z) we get that Hilbert's tenth problem for $\Q$ is undecidable and, in fact, Turing-equivalent to $0'$ by the case of $\Z$, see \cite{DPR, Matiyasevich}. In Section \ref{SecConsequences} we discuss a couple of consequences of the possible existence of a Diophantine definition of $\Z$ in $\Q$.

\subsection{Previous work} The subject started with the first order definition of $\Z$ in $\Q$ given by Julia Robinson in \cite{Robinson}. Several later works approached Hilbert's tenth problem for $\Q$ and a possible Diophantine (i.e. positive existential) definition or interpretation of $\Z$ in $\Q$, see for instance \cite{PoonenDef, Koenigsmann, PoonenSubrings, GFPV}. But the one that is the closest to the current developments is Pheidas's paper \cite{Pheidas}.

The main strategy of using an elliptic curve of rank $1$, local divisibility relations, and an even-order poles predicate $\Phi$ (see Section \ref{SecPhi} for details) is due to Thanases Pheidas who expected that these pieces put together should yield a negative solution to Hilbert's tenth problem over $\Q$, see \cite{Pheidas}  and see also \cite{DemeyervanGeel, CornelissenZahidi} for later work. We call $\Phi$ the \emph{Pheidas predicate} (see Lemma 4.8 in \cite{Pheidas}). 

Our work was motivated by the recent announcement of a negative solution to Hilbert's tenth problem over $\mathbb{Q}$ by OpenAI \cite{OpenAI} which in turn depends on two other of their preprints \cite{OpenAI1, OpenAI2}. In fact, this project started as an attempt to understand that work and we hope that the present paper helps for that purpose. Interestingly, while \cite{OpenAI} does not cite \cite{Pheidas}, the Pheidas predicate in \cite{OpenAI} is called $\Phi$, and we follow this convenient notation in honor of Pheidas (\textgreek{Fe'idas}).

\subsection{About the conjectural inputs} How did OpenAI avoid the use of (X), (Y), and (Z) in \cite{OpenAI}? Let us briefly explain this point. 

Instead of assuming (X), they also released the paper \cite{OpenAI1} on a $2$-converse theorem for elliptic curves. We use (X) via Kummer surfaces while \cite{OpenAI} uses the $2$-converse theorem via Selmer calculations on elliptic curves.

Instead of (Y), they prove a result in the direction of the $abc$ conjecture, but exponentially weaker. This uses as input their companion preprint \cite{OpenAI2} on the Fontaine--Mazur conjecture. Results that turn modularity theorems into height bounds for the $abc$ conjecture go back to the work of Murty--Pasten \cite{MurtyPasten}.

Instead of (Z) they use a version of the Green--Tao theorem \cite{GreenTao} along with a sieving argument. In the present paper we could, in principle, do the same and avoid the use of (Z), but this would add considerable technical complications, too much for our modest purposes.


\section{What if $\Z$ is Diophantine in $\Q$?} \label{SecConsequences}

Let us assume that $\Z$ is Diophantine in $\Q$ in this section. Here are some consequences.

\begin{itemize}

\item[(1)] There is an algebraic variety $X$ over $\Q$ such that the real closure of $X(\Q)$ has infinitely many connected components. Thus, Mazur's conjecture is false \cite{Mazur}.

\item[(2)] $\Q$ has the DPRM property and its Diophantine subsets have uniformly bounded existential rank. See \cite{PastenNotes}.

\item[(3)] There is a Diophantine set $S\subseteq\Q$ whose supremum in $\R$ is transcendental. Thus, Pasten's algebraicity conjecture is false \cite{PastenNotes}.

\item[(4)] The set $S=\{2^n : n\in \N\}$ is Diophantine in $\Q$. Thus, there is a geometrically integral variety $X$ with a morphism  $f:X\to \Pro^1$ over $\Q$ such that $f(X(\Q))$ is an infinite subset of $S$. Taking ramified Kummer covers of a projective closure of $X$ one sees that it does not have the weak Hilbert property of Corvaja and Zannier \cite{CorvajaZannier}, which disproves one of their conjectures there.

\item[(5)] Let ${\bf Exists},\, {\bf Closure},\, {\bf Density}$ be the problems of deciding existence of rational points, computing the Zariski closure of the set of rational points, and deciding Zariski density of rational points, for algebraic varieties over $\Q$. Then, unconditionally, one has
$$
{\bf Closure}\equiv_T {\bf Exists}\oplus{\bf Density}.
$$
Assuming that $\Z$ is Diophantine in $\Q$, and using a uniform Diophantine enumeration of listable sets, ${\bf Closure}$ decides infinitude of listable sets, which is $\equiv_T 0''$. We deduce that both ${\bf Closure}$ and ${\bf Density}$ are undecidable and, furthermore, they are not $\le_T 0'$.
\end{itemize}


\section{The elliptic curve}\label{SecElliptic}

Let $\theta=\sqrt{2}$ and $F=\Q(\theta)$. Define the elliptic curve 
$$
E : \quad y^2 = x^3-(3-\theta)^2x.
$$
Up to a change of variables, this is the same elliptic curve from \cite{OpenAI}.

A direct computer-assisted calculation gives

\begin{lemma} The Mordell--Weil rank of $E(F)$ is $1$ and it has full $F$-rational $2$-torsion. Moreover, there is $P\in E(F)$ such that $2E(F)=\Z\cdot P$.
\end{lemma}
Although $E$ is defined over $F$, it can be expressed with its group structure over $\Q$ in a Diophantine way by separating the rational and irrational parts of all numbers and variables that appear.

The rational function $\xi = -x/y$ is the standard local parameter at the infinity point $O$ of $E$.  The places of bad reduction divide $2$ and $7$. For other primes $w$ of $F$ above a prime $q\ne 2,7$ (normalized by $v_w(q)=1$), for $m,n$ with $mP,nP$ reducing to $O$ at $w$ we have
\begin{equation}\label{EqnEcongr1}
v_w\left(\frac{\xi(nP)}{\xi(mP)}\right) = v_q(n/m)
\end{equation}
and
\begin{equation}\label{EqnEcongr2}
\frac{\xi(nP)}{\xi(mP)} \equiv n/m \bmod w \quad \mbox{if }v_q(n/m)\ge 0
\end{equation}
which can be seen using the formal logarithm of $E$ at $w$.

If $u=a/b$ is a rational number with $a,b$ coprime integers, its logarithmic height is
$$
h(u)=\log \max\{|a|,|b|\}.
$$
The height extends in a standard way to number fields. On $E(F)$ the canonical height $\hat{h}$ is normalized so that as $n$ varies we have
$$
h(x(nP)) = 2n^2\hat{h}(P) + O_E(1).
$$
The value $\hat{h}(P)$ is strictly positive.

Let us say that a prime number $r$ is \emph{tailored} if it satisfies:
$$
r\equiv 1\bmod 8,\quad r\equiv 3\bmod 7,\quad \mbox{ and } (r+1)/2\mbox{ is prime.}
$$
These primes, in particular, split in $F$ and are not divisible by primes of bad reduction of $E$.
\begin{lemma}\label{LemmaTailored} Let $r$ be a tailored prime and let $w,w'$ be the places of $F$ above $r$. Let $d$ and $d'$ be the orders of $E(\F_w)$ and  $E(\F_{w'})$. Then
$$
d+d' = 2(r+1),\quad \gcd(d,d')=4,\quad \mbox{ and }r\nmid dd'.
$$
\end{lemma}
\begin{proof} The reductions of $E$ at $w$ and $w'$ differ by a quadratic twist by $7$. By quadratic reciprocity $(7/r)=(r/7)=(3/7)=-1$ and the first condition follows.

Both $d$ and $d'$ are divisible by $4$ because there is full $F$-rational $2$-torsion and $d/4+d'/4 = (r+1)/2$, which is prime by assumption, and the gcd assertion follows. The last property follows from Hasse's bound and divisibility by $4$.
\end{proof}


\section{Application of (X): Pheidas's predicate}\label{SecPhi}

For $i=1,2$ define $f_{c,i}(X) = c^2X^4+X+i$. Pheidas's unary predicate $\Phi$ in $\Q$ is defined as follows:

We declare that $\Phi(0)$ holds in $\Q$. For $b\ne 0$ we define $\Phi(b)$ as
$$
\exists \alpha,t,c,x,y,z:\quad \alpha=\pm 2b\quad \&\quad c=\alpha+\alpha^3t^4\quad\&\quad cz^2 = f_{c,1}(x) f_{c,2}(y) \quad\&\quad cz\ne 0.
$$
Note that $\Phi$ is Diophantine in $\Q$.
\begin{lemma}\label{LemmaPhi} Pheidas's predicate satisfies the following:
\begin{itemize}
\item[(i)] If $\Phi(b)$ holds, then for every odd prime $q$ with $v_q(b)<0$ (poles) we have that $v_q(b)$ is even.
\item[(ii)] If (X) holds and if $b\in \Z$, then $\Phi(b)$ holds.
\end{itemize}
\end{lemma}
\begin{proof}
Item (i) is a straightforward verification that we leave to the reader.

Under the assumption that (X) holds, item (ii) follows from Theorem A in \cite{HarpazSkorobogatov} applied to the surface
$$
S:\quad Z^2 = f_{c,1}(x)\cdot c f_{c,2}(y)
$$
Under $Z=cz$, this is equivalent to the equation in the definition of $\Phi$. For $c\ne 0$, the surface $S$ is an affine model of a Kummer surface except at points with $z=0$. Under assumption (X), we will show that when $b\ne 0$ is an integer, there are suitable integer values of $\alpha$ and $t$ such that with $c=\alpha+\alpha^3t^4$ the surface $S$ has a Zariski dense set of rational points, which in particular shows that $\Q\models \Phi(b)$.

Under assumption (X), Theorem A from \cite{HarpazSkorobogatov} proves a local-global principle yielding a Zariski dense set of rational points, provided that the following conditions are met:
\begin{itemize}
\item[(i)] $S$ is everywhere locally soluble with $Z\ne 0$.
\item[(ii)] The polynomials $f_{c,1}(x)$ and $f_{c,2}(y)$ are irreducible over $\Q$ with Galois group $S_4$.
\item[(iii)] There are primes $p_1,p_2$ not dividing $c$, such that $v_{p_i}(\Delta_j)=\delta_{ij}$ where $\Delta_j$ is the discriminant of $f_{c,j}$.
\end{itemize}

Assume that $b\in\Z$ is not zero and choose $\alpha=2|b|$. For the moment assume $t$ is an integer; we will later find a suitable value for it.

Everywhere local solubility is a standard application of Hensel's lemma.

Treat $t$ as a variable $T$ and consider $c=c(T)$. For the discriminants, one notes that if $D_i=256i^3c^2-27$ then
$$
\Delta_i = c^4 D_i\quad \mbox{ and }\quad D_2-8D_1=189.
$$
Then one can find $p_1,p_2>7$ where $D_i(T)$ has a simple root modulo $p_i$ and lift it modulo $p_i^2$ to a value that is not a root. This produces an arithmetic progression of integer $t$'s with the required primes.

The polynomials $f_{c(T),i}$ over $\Q(T)$ have the desired Galois group (using the concrete formula for $c$ in terms of $t$) and a Hilbert irreducibility argument over integer arithmetic progressions produces many integer values of $t$ with the same property.

Thus, Theorem A from \cite{HarpazSkorobogatov} applies.
\end{proof}

\section{Application of (Y): multiplicities}

Let us recall the following result of Granville \cite{Granville}.

\begin{theorem} Let $G(x,y)$ be a homogeneous binary form with coefficients in $\Z$ and no repeated factor, of degree $r\ge 3$. Assuming the $abc$ conjecture, for every $\epsilon>0$ one has
$$
\max\{|a|,|b|\}^{r-2-\epsilon}\ll_{G,\epsilon} \rad\, G(a,b)
$$
where $(a,b)$ varies over pairs of coprime integers that are not zeros of $G$.
\end{theorem}

Let $S$ be a finite set of primes containing $2$ and for a rational number $u\ne 0,1,2,3$ define $N_S(u)$ as the product of the primes $q\notin S$ with $v_q(u-i)$ odd for at least one $i=0,1,2,3$. For $u=0,1,2,3$ define $N_S(u)=1$. 

\begin{lemma} \label{Lemmaabc} Assume (Y). Then for every rational number $u$ with $\Q\models \Phi(u)$ we have
$$
h(u)\le 2\log N_S(u) + O_S(1).
$$
\end{lemma}
\begin{proof}
This follows by applying Granville's theorem to
$$
G(x,y) = xy(x-y)(x-2y)(x-3y)
$$
 Indeed, one replaces $x$ and $y$ by the numerator and denominator $a,b$ of $u$ (in reduced form) and we recall that $b$ will only have even exponents outside $S$. The crucial observation is
$$
 \rad\, G(a,b)^2 \ll_{S} |G(a,b)| N_S(u) \ll \max\{|a|,|b|\}^5N_S(u).
$$
Then, under (Y) we have
$$
\max\{|a|,|b|\}^{6-\epsilon}\ll_{S,\epsilon }  \max\{|a|,|b|\}^5 N_S(u)
$$
and taking $\epsilon=1/2$ we get the result.
\end{proof}

\section{Application of (Z): prime representations}

Recall the definition of tailored primes from Section \ref{SecElliptic}.

\begin{lemma}\label{LemmaPrimes} Assume (Z). There is a finite set $K$ of positive integers such that every $A\in \Z$ admits a decomposition
$$
A=u_1+u_2+u_3
$$
with $u_i -j= \eta_i k_{ij} r_{ij}$ for $i=1,2,3$ and $j=0,1,2,3$, where $u_i$ are integers, $(\eta_1,\eta_2,\eta_3)=(1,1,-1)$, $k_{ij}\in K$, and the $r_{ij}$ are tailored primes that can be taken as large as we want.
\end{lemma}
\begin{proof} We only sketch the argument, as it is mostly routine. Consider the degree $1$ polynomials
$$
L_{1j} = x-j,\quad L_{2j} = y-j, \quad L_{3j} =x+y-A+j.
$$
Ideally one would like that for each $j=0,1,2,3$ these polynomials take tailored prime values, which in particular implies that for each prime value $r$ the number $(r+1)/2$ is also prime. This sets the condition that 24 linear polynomials take simultaneous prime values in prescribed congruence classes. This is of course not possible, due to congruence obstructions to small moduli. But after suitable changes of variables and restricting to arithmetic progressions, the new degree $1$ polynomials (now in $1$ variable) no longer have those obstructions and (Z) affords the required primes. Importantly, in this way the finite set $K$ is independent of $A$.
\end{proof}


\section{Local predicates}

This section is a modern version of Pheidas's item b' from \cite{Pheidas}. We later learned that \cite{Farhangi} arrived at a similar computation, and we incorporate here some of those ideas to simplify the argument.

Fix a finite set $S$ containing $2,7$, and all prime divisors of elements of $K$. For $r\in\Q^\times$ define
$$
A_r=\{a\in\Q:\Phi(ra^2)\}\quad\mbox{and}\quad M_r=A_r+\theta A_r\subset F.
$$
Call $(r,\beta,\kappa)\in\Q^3$ a \emph{contact datum} if
\begin{equation}\label{EqnContact}
r\ne0,\quad\beta^2-2=r\kappa,\quad\mbox{ and }\kappa\in A_r.
\end{equation}
Note that $\kappa$ can be expressed in terms of $r$ and $\beta$, but we keep it in the notation for convenience in the computations.

Given a prime $q$, write $w_\beta$ for the place above $q$ selected by $\theta\equiv\beta$; the notation does not reflect dependence on $q$ but this should not lead to confusion. 

For $Y=y_0+y_1\theta$ with $y_i\in\Q$ write $Y_\beta=y_0+\beta y_1$. For $\gamma\in F$ define
\begin{equation}\label{EqnIdeal}
\Ifrak_{r,\beta}(\gamma):\quad \exists X,Y\in M_r, \mu\in A_r:\quad
\gamma Y=rX\quad\&\quad Y\ne0\quad\&\quad Y_\beta+r\mu=1.
\end{equation}

\begin{lemma}\label{LemmaLocal}
Let $(r,\beta,\kappa)$ be a contact datum and let $q\notin S$ be a prime with $v_q(r)>0$ odd. Then $q$ splits in $F$, the number $\beta$ selects a place $w_\beta$ by $\theta\equiv\beta\pmod{w_\beta}$, and
$$
\Q\models \Ifrak_{r,\beta}(\gamma)\quad\Longrightarrow\quad v_{w_\beta}(\gamma)>0.
$$
Conversely, under (X), for any given split prime $q\notin S$ and $\beta\in\Z$ with $\beta^2\equiv2\pmod q$, the triple $(q,\beta,(\beta^2-2)/q)$ is a contact datum and
$$
\Q\models\Ifrak_{q,\beta}(\gamma)\quad\Longleftrightarrow\quad v_{w_\beta}(\gamma)>0.
$$
\end{lemma}

\begin{proof}
Let $e=v_q(r)>0$ be odd. Lemma \ref{LemmaPhi} gives $v_q(a)>-e/2$ for $a\in A_r$, hence $v_w(Z)>-e/2$ for $Z\in M_r$ and $w\mid q$. The contact datum gives $v_q(\beta^2-2)>e/2$, so $q$ splits and writing $v=v_{w_\beta}$ we get
$$
v(\theta-\beta)=v_q(\beta^2-2)>e/2.
$$
If $\Q\models\Ifrak_{r,\beta}(\gamma)$ we take $X,Y,\mu$ as in \eqref{EqnIdeal}. Write $Y=y_0+y_1\theta$ and observe that
$$
Y=1-r\mu+y_1(\theta-\beta)\equiv1\bmod{w_\beta}.
$$
 Thus $Y$ is a unit at $w_\beta$, and $\gamma Y=rX$ yields
$$
v(\gamma)=v(rX)>e/2>0.
$$

Conversely, (X) and Lemma \ref{LemmaPhi} give $\Z\subset A_q$ and $O_F\subset M_q$, so the provided triple is a contact datum and in particular the claimed forward implication holds by the first part. Let us check the converse implication.

Assume $v_{w_\beta}(\gamma)>0$. Then $\gamma/q$ belongs to the local ring at $w_\beta$. The ring $O_F$ is Euclidean and we can find coprime $x,y\in O_F$ such that $x/y = \gamma/q$ is in the local ring at $w_\beta$. Hence $ w_\beta$ does not divide $y$ and since $(q,\beta,(\beta^2-2)/q)$ is a contact datum we conclude that $q\nmid y_\beta$.

Bezout yields $b,\mu\in\Z$ with $by_\beta+\mu q=1$. Taking $X=bx$ and $Y=by$ we see that  $X,Y,\mu$ make \eqref{EqnIdeal} hold.
\end{proof}

Let us mention a key consequence of the previous lemma. For $T\in 2E(F)=\Z\cdot P$ put
$$
\Rfrak_{r,\beta}(T):\quad T=O\, \vee\, \left(T\ne O\, \wedge\, \Ifrak_{r,\beta}(1/x(T))\right).
$$
For $T\ne O$ at a good place $w$, reduction to $O$ is equivalent to $v_w(x(T))<0$.

Thus, for a contact datum, $\Rfrak_{r,\beta}(T)$ forces reduction to $O$ at $w_\beta$ for each $q\notin S$ prime with $v_q(r)>0$ odd (here, $w_\beta$ is the prime above $q$ determined by the contact datum). Under (X), if $r=q\notin S$ is a split prime and $\beta$ is an integer with $\beta^2\equiv 2\bmod q$, then Lemma \ref{LemmaLocal} shows that $\Rfrak_{r,\beta}(T)$ is equivalent to that reduction. 

\section{Defining the integers}

For $Q\in 2E(F)=\Z P$ and $t\in\Q$, let $\Lfrak_{r,\beta}(Q,t)$ assert the existence of $H,J\in 2E(F)\setminus\{O\}$ and $U\in 2E(F)$ such that, with $V_0=H-[4]P$ and $V_1=U-[4]Q$,
\begin{equation}\label{EqnRels}
\begin{aligned}
\Rfrak_{r,\beta}(H),&\quad \Rfrak_{r,\beta}(U), &\quad \Ifrak_{r,\beta}\left(\frac{\xi(U)}{\xi(H)}-t\right),&\\
\Rfrak_{r,-\beta}(J),& \quad\Rfrak_{r,-\beta}(V_i),&\quad \Ifrak_{r,-\beta}\left(\frac{\xi(V_i)}{\xi(J)}\right)&\quad(i=0,1).
\end{aligned}
\end{equation}
The denominators are nonzero because $H,J\ne O$.

\begin{lemma}\label{LemmaLinear}
For a contact datum $(r,\beta,\kappa)$, if $\Lfrak_{r,\beta}(nP,t)$ holds, and if $q\notin S$ has $v_q(r)>0$ odd, then $v_q(t-n)>0$. Conversely, under (X), every tailored prime $r\notin S$ admits integers $\beta,\kappa$ forming a contact datum for which $\Lfrak_{r,\beta}(nP,n)$ holds for all $n\in\Z$.
\end{lemma}

\begin{proof}
Write $H=aP$, $J=bP$, $U=dP$. The triple $(r,-\beta,\kappa)$ is also a contact datum, selecting the conjugate place. The second line of \eqref{EqnRels}, Lemma \ref{LemmaLocal}, and \eqref{EqnEcongr1} give
$$
v_q((a-4)/b)>0,\qquad v_q((d-4n)/b)>0.
$$
Since $b\in\Z\setminus\{0\}$, we have $a\equiv4$ and $d\equiv4n\pmod q$. Thus $q\nmid a$, and \eqref{EqnEcongr2} gives $\xi(U)/\xi(H)\equiv d/a\equiv n\pmod{w_\beta}$. The first line of \eqref{EqnRels} now implies $v_q(t-n)>0$.

Conversely, choose $\beta\in\Z$ with $\beta^2\equiv2\pmod r$ and put $\kappa=(\beta^2-2)/r$. By Lemma \ref{LemmaLocal} this is a contact datum. Let $d_\pm$ be the order of the group of mod $w_{\pm\beta}$ points of $E$. By Lemma \ref{LemmaTailored} we can choose $a\ne0$ with
$$
a\equiv0\pmod{d_+},\qquad a\equiv4\pmod{rd_-}.
$$
Take $H=a\cdot P$, $J=d_-\cdot P$, $U=an\cdot P$. These points and $V_0=(a-4)P$, $V_1=(a-4)nP$ have the required reductions. Since $r\nmid d_-$, Equation \eqref{EqnEcongr1} gives positive valuation for $\xi(V_0)/\xi(J)$ and $\xi(V_1)/\xi(J)$ at $w_{-\beta}$ while \eqref{EqnEcongr2} gives $\xi(U)/\xi(H)\equiv n\pmod{w_\beta}$. Lemma \ref{LemmaLocal} proves \eqref{EqnRels}. This includes zero numerators, since $\xi(O)=0$ and $\Ifrak_{r,\pm\beta}(0)$ holds.
\end{proof}

For $Q\in 2E(F)$ and $t,A\in\Q$, let $\Cfrak(Q,t;A)$ be the (Diophantine) relation defined by the condition that there exist $u_1,u_2,u_3\in \Q$, a choice of $k_{\ell j}\in K$ (with the same $K$  from Lemma \ref{LemmaPrimes}), and $r_{\ell j}, \beta_{\ell j}, \kappa_{\ell j}\in\Q$ for $\ell=1,2,3$ and $j=0,1,2,3$, such that:
\begin{itemize}
\item $A=u_1+u_2+u_3$
\item $\Phi(u_\ell)$ for each $\ell=1,2,3$
\item $u_\ell -j = \eta_\ell k_{\ell j} r_{\ell j}$ for each $\ell=1,2,3$ and $j=0,1,2,3$, where $(\eta_1,\eta_2,\eta_3)=(1,1,-1)$
\item The triples $(r_{\ell j},\beta_{\ell j},\kappa_{\ell j})$ are contact data
\item $\Lfrak_{r_{\ell j},\beta_{\ell j}}(Q,t)$.
\end{itemize}
Set $x_*(Q)=x(Q)$ for $Q\ne O$ and $x_*(O)=0$. Define the Diophantine predicate $\Psi(t)$ by the conditions
\begin{equation}\label{EqnPsi}
\begin{aligned}
\exists Q\in 2E(F),\, \zeta_0,\zeta_1,\zeta_2\in\Q:\quad
\zeta_2\ne0\quad\& \quad\zeta_2x_*(Q)=\zeta_0+\theta\zeta_1\quad\& \quad \Cfrak(Q,t;t^7) \quad\& \quad \bigwedge_{i=0}^2\Cfrak(Q,t;\zeta_i).
\end{aligned}
\end{equation}

\begin{proposition}\label{PropFinal}
Under (X)--(Z), $\Psi(n)$ holds for all $n\in\Z$, and $\Xcal=\{t\in\Q\setminus\Z:\Psi(t)\}$ is finite.
\end{proposition}
\begin{proof} From now on we assume (X)--(Z).

For $t=n\in\Z$, take $Q=nP$ and integral $\zeta_0,\zeta_1,\zeta_2$ with $\zeta_2\ne0$ and $\zeta_2x_*(Q)=\zeta_0+\theta\zeta_1$, using $(0,0,1)$ if $Q=O$. Lemma \ref{LemmaPrimes} decomposes the four integers $\zeta_0,\zeta_1,\zeta_2,n^7$ as in the definition of $\Cfrak(Q,n;n^7)$ and $\Cfrak(Q,n;\zeta_i)$ using tailored primes outside $S$; Lemmas \ref{LemmaPhi} and \ref{LemmaLinear} supply the remaining conditions.

Suppose $\Psi(t)$ holds with $t\notin\Z$ for a choice $Q=nP$ and $\zeta_i\in \Q$. Write $t=A/B$ in lowest terms with $B>0$, and put
$$
\delta=A-nB \ne 0
$$
 For any summand $u=u_\ell$ in the definition of the $\Cfrak$'s in $\Psi(t)$, a prime $q\mid N_S(u)$ has odd positive valuation in some $u-j$, hence in the corresponding $r_{\ell j}$, since the prime divisors of $k_{\ell j}$ lie in $S$. Lemma \ref{LemmaLinear} gives $v_q(t-n)>0$. As $\gcd(\delta,B)=1$, this implies $q\mid\delta$. Thus $N_S(u)\mid\delta$, and Lemma \ref{Lemmaabc} gives $h(u)\le2\log|\delta|+O_{K,S}(1)$. Summing three such terms yields
\begin{equation}\label{EqnBounds}
h(\zeta_i)\le6\log|\delta|+O_{K,S}(1)\quad(0\le i\le2)\quad\mbox{ and }\quad 7h(t)=h(t^7)\le6\log|\delta|+O_{K,S}(1).
\end{equation}
Since $|\delta|\le e^{h(t)}(1+|n|)$, the last inequality implies
\begin{equation}\label{Eqnhtbd}
h(t)\le6\log(1+|n|)+O_{K,S}(1).
\end{equation}
For $n\ne0$, the properties of the canonical height on $E$ together with \eqref{EqnPsi} and \eqref{EqnBounds} give
$$
\begin{aligned}
2\hat h(P)n^2&\le h(\zeta_0)+h(\zeta_1)+h(\zeta_2)+O(1)\le18\log|\delta|+O_{K,S}(1)\\
&\le18h(t)+18\log(1+|n|)+O_{K,S}(1)\le126\log(1+|n|)+O_{K,S}(1).
\end{aligned}
$$
Hence $|n|$ is bounded. Equation \eqref{Eqnhtbd} then bounds $h(t)$, also when $n=0$, leaving finitely many possible $t$.
\end{proof}

Theorem \ref{ThmMain} is now immediate from Proposition \ref{PropFinal}.


\section{Acknowledgements and use of AI}

A.K. was supported by ANID Fondecyt Iniciaci\'{o}n grant 11250060 and  National Center for Artificial Intelligence CENIA FB210017. H.P. was supported by ANID Fondecyt Regular grant 1230507 from Chile. X.V. was supported by Proyecto 2025001477INV from Universidad de Concepci\'on.

This project started as an effort to understand OpenAI's AI-generated preprint \cite{OpenAI}, by using frontier AI models. We had the idea of replacing the difficult arithmetic inputs by (X), (Y), and (Z), as well as translating the original argument of \cite{OpenAI} into a finitary version; the details of the computation were assisted by AI. We first obtained a conditional proof of a Diophantine \emph{model} of $\Z$ in $\Q$ under (X), (Y), and (Z), but then on November 5 of 2026 both AI-generated preprints \cite{Farhangi} and \cite{Glazer} appeared, proposing a Diophantine definition of $\Z$ in $\Q$ (the first assumes correctness of three OpenAI papers, the second is self-contained). This made us reconsider our Diophantine model and a small modification yielded the current Diophantine definition of $\Z$ in $\Q$, conditional on (X), (Y), (Z). The core strategy of using an elliptic curve of rank $1$, local divisibility relations, and the even-order poles predicate is due to Thanases Pheidas \cite{Pheidas}. Thus, this note owes an intellectual debt to the work of Pheidas, OpenAI, Farhangi, and Glazer.


\end{document}